\documentclass{amsart}
\usepackage{amsmath,amssymb,amsxtra,amsthm,amscd}
\usepackage{lmodern}
\usepackage[all,cmtip]{xy}
\usepackage{microtype}
\usepackage{todonotes}
\usepackage{multirow}

\usepackage{color}

\newif\ifShowLabels
\ShowLabelstrue
\newcommand{\TeXref}[1]{
\marginpar{\scriptsize \texttt{#1}}}

\DeclareMathOperator{\asdim}{asdim}

\DeclareMathOperator{\diam}{diam}

\DeclareMathOperator{\Out}{Out}

\DeclareMathOperator{\id}{id}

\DeclareMathOperator*{\one}{1}
\newcommand{\onehatplace}[1]
{ \one^{\substack{#1 \\ \frown}} }

\DeclareMathOperator*{\bones}{\times}
\newcommand{\undertimes}[1]
{ \bones_{#1} }

\DeclareMathOperator*{\bowl}{\cup}
\newcommand{\undercup}[1]
{ \bowl_{#1} }

\DeclareMathOperator*{\arch}{\cap}
\newcommand{\undercap}[1]
{ \arch_{#1} }

\newcommand{\pull}
{\!\!\! -\!\!\! -\!\!\! -\!\!\!}

\DeclareMathOperator*{\holimprep}{holim}                       
\newcommand{\holim}[1]%
{\displaystyle\holimprep_{\substack{\leftarrow \pull - \\ #1}} \, }

\DeclareMathOperator*{\hocolimprep}{hocolim}                   
\newcommand{\hocolim}[1]%
{\displaystyle\hocolimprep_{\substack{- \pull \rightarrow \\ #1}} \, }

\DeclareMathOperator*{\plainlim}{lim}                           
\newcommand{\contralim}[1]%
{\displaystyle\plainlim_{\substack{\leftarrow \pull - \\ #1}} \, }

\DeclareMathOperator*{\plaincolim}{colim}                       
\newcommand{\colim}[1]%
{\displaystyle\plaincolim_{\substack{- \pull \rightarrow \\ #1}} \, }

\DeclareMathOperator*{\laxlimplain}{laxlim}                     
\newcommand{\laxlim}[1]%
{\displaystyle\laxlimplain_{\substack{\leftarrow \pull - \\ #1}} \, }

\providecommand{\bysame}{\makebox[3em]{\hrulefill}\thinspace}

\theoremstyle{plain}
\newtheorem{Thm}{Theorem}[section]

\newtheorem*{MainThm}{Main Theorem}

\newtheorem{Cor}[Thm]{Corollary}

\newtheorem{Lem}[Thm]{Lemma}
\newtheorem{Prop}[Thm]{Proposition}

\theoremstyle{definition}
\newtheorem{Def}[Thm]{Definition}
\newtheorem{Q}[Thm]{Question}

\newtheorem{Ex}[Thm]{Example}

\newtheorem*{Exs}{Examples}

\newtheorem{Rem}[Thm]{Remark}

\theoremstyle{remark}
\newtheorem{Not}[Thm]{Notation}

\newtheoremstyle{freestylethm}{6pt}{6pt}{\itshape}{}%
                {\bfseries}{}{.5em}{\thmnote{#3}}
\theoremstyle{freestylethm}

\newcommand{\SecRef}[2]{\section{#1}\label{S:#2}%
\ifShowLabels \TeXref{{S:#2}} \fi}

\newcommand{\refT}[1]{\textup{\ref{T:#1}}}

\newcommand{\refD}[1]{\textup{\ref{D:#1}}}

\newenvironment{ThmRef}[1]%
{ \begin{Thm} \label{T:#1}
\ifShowLabels \TeXref{T:#1} \fi }%
{ \end{Thm} }
\newenvironment{DefRef}[1]%
{ \begin{Def} \label{D:#1}
\ifShowLabels \TeXref{D:#1} \fi }%
{ \end{Def} }
\newenvironment{QRef}[1]%
{ \begin{Q} \label{Q:#1}
\ifShowLabels \TeXref{Q:#1} \fi }%
{ \end{Q} }
\newenvironment{LemRef}[1]%
{ \begin{Lem} \label{L:#1}
\ifShowLabels \TeXref{L:#1} \fi }%
{ \end{Lem} }
\newenvironment{CorRef}[1]%
{ \begin{Cor} \label{C:#1}
\ifShowLabels \TeXref{C:#1} \fi }%
{ \end{Cor} }
\newenvironment{RemRef}[1]%
{ \begin{Rem} \label{R:#1}
\ifShowLabels \TeXref{R:#1} \fi }%
{ \end{Rem} }
\newenvironment{PropRef}[1]%
{ \begin{Prop} \label{P:#1}
\ifShowLabels \TeXref{P:#1} \fi }%
{ \end{Prop} }
\newenvironment{ExRef}[1]%
{ \begin{Ex} \label{E:#1}
\ifShowLabels \TeXref{E:#1} \fi  }%
{ \end{Ex} }
{ \begin{Not} \label{N:#1}
\ifShowLabels \TeXref{N:#1} \fi }%
{ \end{Not} }

{ \begin{Thm} [#2]\label{T:#1}
\ifShowLabels \TeXref{T:#1} \fi }%
{ \end{Thm} }
\newenvironment{DefRefName}[2]%
{ \begin{Def} [#2]\label{D:#1}
\ifShowLabels \TeXref{D:#1} \fi }%
{ \end{Def} }
{ \begin{Lem} [#2]\label{L:#1}
\ifShowLabels \TeXref{L:#1} \fi }%
{ \end{Lem} }
{ \begin{Cor} [#2]\label{C:#1}
\ifShowLabels \TeXref{C:#1} \fi }%
{ \end{Cor} }
{ \begin{Rem} [#2]\label{R:#1}
\ifShowLabels \TeXref{R:#1} \fi }%
{ \end{Rem} }
{ \begin{Prop} [#2]\label{P:#1}
\ifShowLabels \TeXref{P:#1} \fi }%
{ \end{Prop} }
{ \begin{Ex} [#2]\label{E:#1}
\ifShowLabels \TeXref{E:#1} \fi }%
{ \end{Ex} }

 \ShowLabelsfalse

\begin{document}

\title[Extension properties of APC and FDC]
{Extension properties of asymptotic property C and Finite Decomposition Complexity}
\author[Susan Beckhardt]{Susan Beckhardt}
\address{Department of Mathematics and Statistics\\ SUNY\\ Albany\\ NY 12222}
\email{sbeckhardt@albany.edu}
\author[Boris Goldfarb]{Boris Goldfarb}
\address{Department of Mathematics and Statistics\\ SUNY\\ Albany\\ NY 12222}
\email{bgoldfarb@albany.edu}
\date{\today}

\begin{abstract}
We prove extension theorems for several geometric properties such as {asymptotic property C} (APC), {finite decomposition complexity} (FDC), strict finite decomposition complexity (sFDC) which are weakenings of Gromov's finite asymptotic dimension (FAD).  

The context of all theorems is a finitely generated group $G$ with a word metric and a coarse quasi-action on a metric space $X$. 
We assume that the quasi-stabilizers have a property $P_1$, and $X$ has the same or sometimes a weaker property $P_2$.  
Then $G$ also has property $P_2$.  

We show some sample applications, discuss constraints to further generalizations, and illustrate the flexibility that the weak quasi-action assumption allows.
\end{abstract}

\maketitle

\tableofcontents

\SecRef{Introduction}{S}

There are a number of geometric conditions of asymptotic or coarse nature that one may impose on a metric space or a group with a word metric.  They have been exploited to great success in the recent work of many authors on the Novikov and Borel conjectures in manifold topology, the problems which are known to have subtle dependence on the asymptotic properties of unbounded universal covers of compact manifolds.

The properties we focus on are best viewed as weakenings of the celebrated property introduced by Gromov in \cite{mG:93} called finite asymptotic dimension (FAD). The groups with asymptotic property C (APC), finite decomposition complexity (FDC), strict finite decomposition complexity (sFDC), and some related properties such as fsFDC that we review in the next section are intermediate between the very large class FAD and the groups with Yu's property A, as shown in \cite{aDmZ:12}.  The well-known example of a group with infinite asymptotic dimension, the restricted wreath product of two copies of the integers $\mathbb{Z}$, has FDC.
Further examples of FDC groups established in \cite{eGrTgY:12} are
all finitely generated subgroups of $GL_n (k)$, where $k$ is a field.
There are other infinite dimensional groups such as Thompson's group, Grave's group \cite{bG:03}, Gromov's random groups. It is unknown where they or other groups of much interest such as $\textrm{Out} (F_n)$ belong in this hierarchy.

Our main theorem is a generalization of a theorem of Bell and Dranishnikov about groups acting on finite dimensional spaces \cite[Theorem 2]{gBaD:06}.

We refer the reader to the next section for the precise definition of a coarse quasi-action on a metric space.  It is a significant and useful generalization of actions by isometries and, further, quasi-isometries.  Quasi-actions by quasi-isometries in particular are the framework for some fundamental questions in geometric group theory, cf. \cite{cD:09,mK:14,lMmSkW:03}.

\begin{MainThm}
Let $G$ be a finitely generated group with a coarse quasi-action on a metric space $X$.
If $X$ has FAD (respectively, APC) and the quasi-stabilizers have asymptotic dimension bounded from above by some number $n \ge 0$  then $G$ has FAD (respectively, APC) with respect to a word metric.
If $X$ and all quasi-stabilizers of the action have FDC (respectively, sFDC) then $G$ has FDC (respectively, sFDC).\end{MainThm}

The theorem has applications to extensions of groups via standard constructions. 
We discuss these results in the last section.

The authors would like to thank Marco Varisco for discussions related to the subject of this paper and the referees for comments and suggestions.

\SecRef{Background from coarse geometry}{CGPREL}

We give a quick review of some coarse geometric finiteness conditions mostly inspired by the covering dimension in topology.  In the second half of the section, we define coarse quasi-actions on metric spaces and point out their most important occurrances in the literature.

\medskip

Given a subset $S$ of a metric space $X$, we will use the notation $S[b]$ for the $b$-enlargement of $S$, that is the subset $\{ x \in X \mid d(x,s) \le b \ \mathrm{for \ some} \ s \in S \}$.
So, in particular, the metric ball centered at $x$ with radius $r$ is denoted by $x[r]$.

Also, given a number $R > 0$, a collection of disjoint subsets $S_{\alpha}$ of $X$ is called $R$-\textit{disjoint} if $S_{\alpha} [R]$ is disjoint from the union $\bigcup_{\beta \ne \alpha} S_{\beta}$, for all $\alpha$.

\begin{DefRefName}{ASDIM}{FAD}
The \textit{asymptotic dimension} of a metric space $X$ was defined by Gromov as
the smallest number $n$ such that for a given positive number $0 < R$ there are $n + 1$
uniformly bounded $R$-disjoint families of subsets of $X$
which together cover $X$.
If such number exists, one says that $X$ has \textit{finite asymptotic dimension}.
\end{DefRefName}

Asymptotic property C defined by Dranishnikov \cite{aD:00} is a weakening of FAD.

\begin{DefRefName}{APC}{APC}
A metric space $X$ has the \textit{asymptotic property C}
if for every sequence of positive numbers $0 < R_1 \le R_2 \le \ldots$ there exists a natural number $n$ and
uniformly bounded $R_i$-disjoint families $\mathcal{W}_i$, $1 \le i \le n$, such that the union of all $n$ families is a covering of $X$.
\end{DefRefName}

On the other hand, the following is one of the equivalent definitions of the finite decomposition complexity from Guentner--Tessera--Yu \cite{eGrTgY:12}.

Let $\mathcal{X}$ and $\mathcal{Y}$ be two families of metric spaces.  Let $R > 0$.
The family $\mathcal{X}$ is called $R$-\textit{decomposable over} $\mathcal{Y}$ if for any space $X$ in $\mathcal{X}$ there are collections of subsets $\{ U_{1,\alpha} \}_{\alpha \in A}$, $\{ U_{2,\beta} \}_{\beta \in B}$ such that
\[
X = \bigcup_{\stackrel{\scriptstyle i=1,2}{\scriptstyle \gamma = A \cup B}} U_{i,\gamma},
\]
each $U_{i,\gamma}$ is a member of the family $\mathcal{Y}$, and each of the collections $\{ U_{1,\alpha} \}$ and $\{ U_{2,\beta} \}$ is $R$-disjoint.
A family of metric spaces is called \textit{bounded} if there is a uniform bound on the diameters of the spaces in the family.

This definition is in terms of a winning strategy for the following game between two players.
The families of metric spaces that appear in the decompositions are families of metric subspaces of $X$.
In round number 1 the first player selects a number $R_1 > 0$, the second player has to select a family of metric spaces $\mathcal{Y}_1$ and an $R_1$-decomposition of $\{ X \}$ over $\mathcal{Y}_1$.
In each succeeding round number $i$ the first player selects a number $R_i > 0$, the second player has to select a family of metric spaces $\mathcal{Y}_i$ and an $R_i$-decomposition of $\mathcal{Y}_{i-1}$ over $\mathcal{Y}_i$.
The second player wins the game if for some finite value $k$ of $i$ the family $\mathcal{Y}_k$ is bounded.

\begin{DefRefName}{FDC}{FDC}
A metric space $X$ has \textit{finite decomposition complexity} if the second player possesses a winning strategy in every game played over $X$.
\end{DefRefName}

The following property was defined by Dranishnikov and Zarichnyi in \cite{aDmZ:12}.

\begin{DefRefName}{DZP}{sFDC}
A metric space $X$ has \textit{straight finite decomposition complexity} if, for any sequence $R_1 \le R_2 \le \ldots$
of positive numbers, there exists a finite sequence of metric families
$\mathcal{V}_1$, $\mathcal{V}_2$, \ldots, $\mathcal{V}_n$ such that $\{ X \}$ is $R_1$-decomposable over $\mathcal{V}_1$, $\mathcal{V}_1$ is $R_2$-decomposable over $\mathcal{V}_2$, etc., and the family $\mathcal{V}_n$ is bounded.
\end{DefRefName}

It is known that the weakest of the four conditions is sFDC and that all four are stronger 
than Yu's property A.

There are also variations on sFDC in the literature.  For example, Ramras and Ramsey \cite{RR} define certain weakenings called $k$-fold straight finite decomposition complexity and weak straight finite decomposition complexity.

\bigskip

Let $X$ and $Y$ be metric spaces with metric functions $d_X$ and $d_Y$.
We will assume that the metrics are proper, in the sense that closed bounded subsets of $X$ and $Y$ are compact.

\begin{DefRef}{Bornol2}
A map $f \colon X \to Y$
between proper metric spaces is called \textit{bornologous} if there is a real positive non-decreasing function $\ell$
such that
\begin{equation}
d_X (x_1, x_2) \le r \ \Longrightarrow
d_Y (f(x_1), f(x_2)) \le \ell(r).  \notag
\end{equation}
The same kind of map is called \textit{proper} if $f^{-1} (S)$ is a bounded subset of $X$ for
each bounded subset $S$ of $Y$.
We say $f$ is a \textit{coarse map} if it is both bornologous and proper.

We will say that a function $f \colon X \to X$ is \textit{bounded} if
$d_X (x, f(x)) \le D$ for all $x \in X$ and a fixed $D \ge 0$.
A bornologous map $f$ is a \textit{coarse equivalence} if there is a bornologous map $g \colon Y \to X$ such
that $f \circ g$ and $g \circ f$ are bounded maps.
\end{DefRef}

\begin{Exs}
Any bounded function $f \colon X \to X$ is coarse.
In fact, it is a coarse
equivalence using $\ell (r) = r + 2D$ for both $f$ and its coarse inverse.

The isometric embedding of a metric subspace is a coarse map.
An isometry, which is a bijective isometric map, is a coarse equivalence.
An isometric embedding onto a subspace that has the property that its bounded enlargement is the whole target metric space is also a coarse equivalence.

There are many more maps in the literature that are coarse.
For example, proper eventually Lipschitz maps are coarse.
Quasi-isometries are coarse maps.  Quasi-isometries onto commensurable subspaces are coarse equivalences.
\end{Exs}

It is a basic construction in geometric group theory that allows to treat a group $G$ equipped with a finite generating set as a metric space with the word-length metric.
This metric makes $G$ a proper metric space with an action by $G$ via left
multiplication.
Different generating sets give quasi-isometric metrics and so coarsely equivalent metric spaces.

The following well-known basic fact due to Shvarts and Milnor is known as ``Milnor's lemma''.
 
\begin{ThmRef}{ML}
Suppose $X$ is a path metric space and $G$ is a group acting properly and cocompactly
by isometries on $X$.  Then $G$ is coarsely equivalent to $X$.
\end{ThmRef}

The coarse equivalence is given by the map $g \mapsto g x_0$ for any
point $x_0$ of $X$.

For example, if $K$ is a finite complex with the fundamental group $G = \pi_1 (K)$,
the inclusion of any orbit of $G$ in the universal cover of $K$ is a coarse equivalence for any choice
of the generating set of $G$.

\medskip

A \textit{coarse quasi-action} is designed to describe situations where elements of a group act on a metric space via a coarse equivalences.  We should point out that it is quite a bit weaker than the categorical notion of a group action in the coarse category.

\begin{DefRef}{UCA}
A \textit{coarse quasi-action} of a group $G$ on a metric space $X$ is an assignment of a bornologous function $f_g \colon X \to X$ for each element $g \in G$ so that the following conditions are satisfied:
\begin{enumerate}
\item all $f_g$ are bornologous with respect to a uniform choice of the function $\ell$,
\item there is a number $A \ge 0$ such that $d(f_{\id}, \id_X) \le A$ in the sup norm,
\item there is a number $B \ge 0$ such that $d (f_g \circ f_h, f_{gh}) \le B$ in the sup norm for all elements $g$ and $h$ in $G$.
\end{enumerate}
As a particular consequence of (2) and (3), all compositions $f_g \circ f_{g^{-1}}$ are $(A + B)$-close to the identity.  This implies that all $f_g$ are coarse equivalences.
\end{DefRef}

The following facts make clear that the maps defining a coarse quasi-action are indeed coarse maps.

\begin{PropRef}{cece}
	Suppose $f \colon X \to Y$ is a map such that there is a bornologous map $h \colon Y \to X$ so that $h \circ f$ is bounded, then $f$ is proper.
\end{PropRef}

\begin{proof}
	Given a bounded subset $S \subset X$, we have $S \subset y[B]$ for some $y \in Y$ and a number $B \ge 0$.  For an arbitrary $x \in f^{-1} (S)$, $d (f(x),y) \le B$, so $d (h(f(x)), h(y)) \le \ell (B)$ if $h$ is controlled by $\ell$.
	Now $d(x,h(f(x))) \le D$ for some $D \ge 0$ independent from the choice of $x$.  This shows $d(x,g(y)) \le \ell(B) +D$, so $\diam (S) \le 2\ell(B) + 2D$.
\end{proof}

\begin{CorRef}{cor to ecec}
	The following properties are equivalent for a map $f \colon X \to Y$:
	\begin{enumerate}
		\item $f$ is a coarse equivalence,
		\item $f$ is a proper coarse equivalence with a proper coarse inverse.
	\end{enumerate}
\end{CorRef}

Coarse quasi-action is a generalization of the notion of \textit{quasi-action} central to the fundamental problem of quasi-isometry classification of finitely generated groups, cf.\,\cite{cD:09,mK:14,lMmSkW:03}.  Quasi-isometries are coarse equivalences with linear control functions.  Quasi-actions appear most naturally in a partial converse to Milnor's Lemma. 
This is a curious term because it is a ``kind of action'' by ``quasi''-isometries, so both meanings of ``quasi'' get (inconveniently?) conflated.

To see a simple example of a quasi-action, consider a metric subset $C$ embedded in $X$.  Suppose a group $G$ acts by isometries on $C$.  If $C$ is commensurable with $X$, we can extend the action of $G$ on $C$ to a quasi-action on $X$ by the following device.  If $B \ge 0$ is a commensurability constant, let $\phi \colon X \to C$ be any function bounded by $B$.  Then we define $f_g (x)$ to be the composition $g (\phi (x))$.  All resulting maps are coarse equivalences with $\ell (x) = x + 2B$ and  the compositions $f_g \circ f_{g^{-1}}$ all $2B$-close to the identity.
This example illustrates a generalization of a well-known construction from geometric group theory: if a metric space is coarsely equivalent to a finitely generated group $G$ with a word metric then the left multiplication action in $G$ can be quisiconjugated to give a coarse quasi-action on $X$.

Another source of quasi-actions of interesting groups on well-understood geometries is a number of recent constructions of actions of groups on quasi-spaces, that is spaces that are quasi-isometric to familiar geometries.  For example, Bestvina, Bromberg, and Fujiwara \cite{mBkBkF:15} construct actions of the mapping class group and the outer automorphism group of a free group of rank $> 1$ on quasi-trees, which can be translated as quasi-actions on trees. 

The general \textit{coarse quasi-actions} are gaining prominence in applications to algebraic $K$-theory.  In that subject, coarse maps are precisely the natural maps between metric spaces that induce maps of bounded $K$-theory spectra $K(X,R)$ built from free $R$-modules parametrized over the metric space $X$.  So $K(X,R)$ is a functor on the coarse category of metric spaces and coarse maps.

  An important instance of a coarse quasi-action that is not by quasi-isometries is in section 2.2 of \cite{gCbG:13}. 
  In that work an action by (not necessarily bounded) isometries needed to be converted to an action by bounded coarse maps via a change of metric in $X$.  For example, this general construction converts the left multiplication action on any group with a word metric to a bounded coarse quasi-action.  The resulting quasi-action is no longer by isometries or even quasi-isometries unless the group is abelian.

  \medskip

The point of view in this paper is that in the search for new explicit examples of groups with properties APC and FDC the greatest benefit is from the weakest assumptions on the action.

\SecRef{Extension theorems for APC, FDC, sFDC, etc.}{EXTTHMS}

Throughout this section, $G$ is a finitely generated group with a coarse quasi-action on a metric space $X$.

We will use $\ell$ as a control function for all $g$ in $G$.

The first theorem is a consequence of \cite[Theorem 1]{gBaD:06}.  In fact, our proof follows closely that of \cite[Theorem 2]{gBaD:01} with modifications to accommodate weaker assumptions on the action.  We include a short proof here in order to compare to the more elaborate proof of Theorem \refT{Susan}.

Let $x_0$ be a chosen base point in $X$.  An \textit{$R$-quasi-stabilizer} $W_R (x_0)$ of $x_0$ is the subset of those elements $g$ in $G$ with the property $d(g (x_0),x_0) \le R$. 

We define $\lambda = \max \{ d(s(x_0), x_0) \mid s \in S \}$, where $S$ is a finite generating set for $G$.
There is a map $\pi \colon G \to X$ given by $\pi (g) = g(x_0)$.  If the action of $G$ on $X$ is by isometries, this map is $\lambda$-Lipschitz.  

\begin{LemRef}{pi_is_lipschitz}
For a coarse quasi-action of $G$ on $X$, $\pi$ is $(\ell (\lambda) + B)$-Lipschitz.
\end{LemRef}

\begin{proof}
Let $g,h \in G$. We proceed by induction on $d(g,h)$. First suppose $d(g,h) = 1$; therefore $g = hs$ for some $s \in S$. Then
\begin{align}
d(\pi(g),\pi(h)) & = d(f_{hs}(x_0), f_h(x_0)) \notag \\
& \leq d(f_h(f_s(x_0)),f_h(x_0)) + B \notag \\
& \leq \ell (d(f_s(x_0), x_0)) + B \notag \\
& \leq \ell (\lambda) + B. \notag 
\end{align}

Now suppose $d(\pi(g),\pi(g')) \leq (\ell (\lambda) + B)k$ whenever $d(g,g') \leq k$. Suppose $d(g,h) = k+1$, so that $g = hs_1 \dots s_k s_{k+1}, s_i \in S$. Then 
\begin{align}
d(\pi(g),\pi(h)) & = d(\pi(hs_1 \dots s_k s_{k+1}),\pi(h)) \notag \\
& \leq d(\pi(hs_1 \dots s_k s_{k+1}), \pi(hs_1 \dots s_k)) + d(\pi(hs_1 \dots s_k), \pi(h)) \notag \\
& \leq \ell(\lambda)+B + (\ell(\lambda)+B)k \notag \\
& = (\ell(\lambda)+B)(k+1). \notag
\end{align}

This shows that for all $g,h \in G$, $d(g,h) \leq k$ gives $d(\pi(g),\pi(h)) \leq (\ell(\lambda)+B)k$ as desired.
\end{proof}

\begin{ThmRef}{Props}
We assume that for all $R \ge 0$ the quasi-stabilizers of $x_0$ satisfy
$\asdim (W_R (x_0)) \le n$ for some common number $n \ge 0$.
If $X$ has FAD then $G$ has FAD.
\end{ThmRef}

\begin{proof}

The orbit $G x_0$ is a subset of $X$ and so has FAD.  This allows us to assume without loss of generality that the action on $X$ is transitive.

Let $\lambda = \max \{ d(s(x_0), x_0) \mid s \in S \}$, where $S$ is a finite generating set for $G$.
There is a map $\pi \colon G \to X$ given by $\pi (g) = g(x_0)$.  If the action of $G$ on $X$ is by isometries, this map is $\lambda$-Lipschitz.  In our case, $d (\pi(g), \pi (gs)) = d(g(x_0), gs (x_0)) \le \ell (d(x_0, s(x_0)) \le \ell (\lambda)$.
So $\pi$ is $\ell (\lambda)$-Lipschitz.

Suppose we have $\asdim (X) \le k$.
Given any $r > 0$, there are $\ell (\lambda) r$-disjoint, $T$-bounded families $\mathcal{F}_0, \mathcal{F}_1, \ldots, \mathcal{F}_k$ which cover $X$.  Let us consider an element $g$ and $x = g(x_0)$.
Let $F$ be a $T$-bounded subset from one of the covering families $\mathcal{F}_i$ with $x \in F$.
We know that $d (g^{-1} (x), x_0) \le A + B$ (where $A$ and $B$ are the constants from Definition \refD{UCA}), so we get $g^{-1} (F) \subset g^{-1} (x[T]) \subset g^{-1} (x) [\ell(T)] \subset  x_0 [A + B + \ell(T)]$.
Therefore, straight from the definition of $\pi$, $g^{-1} \pi^{-1} (F)$ is contained in $W_{A + 2B + \ell(T)} (x_0) = \pi^{-1} (x_0 [A + 2B + \ell(T)])$.
We will denote this particular quasi-stabilizer simply as $W$.
 
By our assumption, $\asdim (W) \le n$ for some $n \ge 0$, so there are $n+1$ families $\mathcal{A}_0, \mathcal{A}_1, \ldots, \mathcal{A}_n$ which cover $W$, which are $r$-disjoint for the given $r$, and which are uniformly bounded by some $K \ge 0$.

For each $F \in \mathcal{F}_i$, choose an element $g_F \in \pi^{-1} (F)$.
Left multiplication by any element in $G$ is an isometry, so the formula $g_F \left(g_F^{-1} \pi^{-1} (F) \cap \mathcal{A}_j \right)$ gives families $\mathcal{A}_{F,j}$ which cover $\pi^{-1} (F)$, which are $r$-disjoint, and which are $K$-bounded. 
Now it is clear that the covering of $G$ by the $r$-disjoint $K$-bounded families
\[
\mathcal{W}_{i,j} = \{ g_F \left(g_F^{-1} \pi^{-1} (F) \cap {A} \right) \mid F \in \mathcal{F}_i, A \in \mathcal{A}_j \}
\]
shows that $\asdim (G) \le (k+1)(n+1)-1$.
\end{proof}

\begin{RemRef}{REMx}
The assumptions do not explicitly require that there is a uniform bound on the asymptotic dimension of all quasi-stabilizers of all $x$ in $X$.  Instead, a geometric property of the action pointed out in the middle of the proof allows to isometrically embed all pullbacks
$\pi^{-1} (F)$ of subsets $F$ with diameter bounded by $T$ in the $( A+2B+\ell(T) )$-quasi-stabilizer $W$ of $x_0$, which is treated as a common ``chopping block''.  This says in particular that all subsets $\pi^{-1} (F)$ of $G$ have asymptotic dimension bounded from above by the same number $n$.

We want to insert a technical comment here.  Unlike the situation with the action of $G$ by isometries, it cannot be assumed that for every subset $F$ with diameter bounded by $T$ there is a number $K$ so that $F \subset g (x_0 [K])$ for some group element $g$.
We were only able to guarantee that $g^{-1} (F) \subset x_0 [K]$ for $K=A+2B+\ell(T)$.
\end{RemRef}

\begin{ThmRef}{Susan}
If for all quasi-stabilizers the asymptotic dimension is uniformly bounded by $n\ge 0$, and $X$ has APC, then $G$ has APC.
\end{ThmRef}

\begin{proof}
Some features of the proof of Theorem \refT{Props} should be borrowed without change.  So we have the $\ell (\lambda)$-Lipshitz projection $\pi \colon G \to X$.  This time the orbit $G x_0$ inherits APC, so we can assume $\pi$ is onto.

Let $0 < r_0 < r_1 < r_2 < \dots$ be a sequence of real numbers.  This allows to generate a new sequence 
$0 < \ell(\lambda) r_{n+1} < \ell(\lambda) r_{2(n+1)} < \ell(\lambda) r_{3(n+1)} < \dots$
Since $X$ has APC, we can choose finitely many uniformly $T$-bounded families $\mathcal{F}_0, \mathcal{F}_1, \dots, \mathcal{F}_m$ which cover $X$ and where each $\mathcal{F}_i$ is $\ell(\lambda) r_{(i+1)(n+1)}$-disjoint for $0 \leq i \leq m$.
Now the pulled-back families $\pi^{-1}(\mathcal{F}_i) = \{ \pi^{-1}(F) \mid F \in \mathcal{F}_i \}$ cover $G$, and each family $\pi^{-1}(\mathcal{F}_i)$ is $r_{(i+1)(n+1)}$-disjoint, though in general their elements are not bounded as subsets of $G$. 

Given a subset $F \in \mathcal{F}_i$ and an element $g_F$ such that $g_F (x_0) \in F$, we have seen that 
$g_F^{-1} (F) \subset g_F^{-1} (x) [\ell(T)] \subset  x_0 [A + B + \ell(T)]$.
Therefore, $g_F^{-1} \pi^{-1} (F)$ is contained in the quasi-stabilizer $W = W_{A + 2B + \ell(T)} (x_0)$.
By the special assumption, $\asdim (W) \le n$, so there are $n+1$ families $\mathcal{A}_0, \mathcal{A}_1, \ldots, \mathcal{A}_n$ which cover $W$, which are $r_{(m+1)(n+1)}$-disjoint and which are uniformly bounded.

We cover $G$ by $(m+1)(n+1)-1$ families of subsets $\mathcal{W}_k$ as follows: for every $0 \leq i \leq m$, $0 \leq j \leq n$, let 
\[
\mathcal{W}_{i(n+1)+j} = \{ g_F \left(g_F^{-1} \pi^{-1} (F) \cap {A} \right) \mid F \in \mathcal{F}_i, A \in \mathcal{A}_j \}.
\]

Notice that
\begin{enumerate}
\item{For each $k$, $0 \leq k \leq (m+1)(n+1)-1$, there is exactly one pair $(i,j)$ such that $k=i(n+1)+j$.}
\item{$\bigcup_{k} \mathcal{W}_k$ is a cover of $G$ because  $\bigcup_{i} \pi^{-1}(\mathcal{F}_i)$ is a cover of $G$.}
\item{Each $\mathcal{W}_{i(n+1)+j}$ is uniformly bounded because its elements are subsets of isometric translations of elements of $\mathcal{A}_j$, which are uniformly bounded.}
\item{Each $\mathcal{W}_{i(n+1)+j}$ is $r_{i(n+1)+j}$-disjoint by the following argument. 
Let 
\[
g_F \left(g_F^{-1} \pi^{-1} (F) \cap {A} \right) \neq g_{F'} \left(g_{F'}^{-1} \pi^{-1} ({F'}) \cap {A} \right) \in \mathcal{W}_{i(n+1)+j}.
\] 
If $F \neq F'$, then $d_G(\pi^{-1}(F), \pi^{-1}(F') \geq r_{(i+1)(n+1)} \geq r_{i(n+1)+j}$ since the family $\pi^{-1}(\mathcal{F}_i)$ is $r_{(i+1)(n+1)}$-disjoint. Otherwise, if $F = F'$ but $A \neq A'$, then $d_G(g_FA, g_FA') \geq r_{(m+1)(n+1)} \geq r_{i(n+1)+j}$ because the family $\mathcal{A}_j$ is $r_{(m+1)(n+1)}$-disjoint.}
\end{enumerate}

Therefore $G$ has APC, as desired.
\end{proof}

It is remarkable that the extension theorems for FDC and sFDC have fewer geometric demands and are easier than the case of APC.  In the case of the action by isometries, a proof for FDC was given by Guentner, Tessera, and Yu \cite{eGrTgY:13}. A proof for sFDC was given by Bell and Moran \cite{gBdM:15}.

\begin{ThmRef}{OUR}
If $X$ has FDC (respectively, sFDC) then $G$ has FDC (respectively, sFDC).
\end{ThmRef}

\begin{proof}
We start with FDC.  
Since $X$ has FDC, then there is always a winning strategy that allows the second player to win the game on $X$ with a bounded family $\mathcal{Y}_k$ in a certain number of steps $k$.  Recall that $\pi \colon G \to X$ is $\ell (\lambda)$-Lipschitz. Every time player one calls out a number $R_i$ for $1 \le i \le k$, player two computes $\overline{R_i} = \max \{ 1, \ell (\lambda) \} R_i$ and uses this number as data for a winning strategy over $X$.  At every step, player two returns the following family as the response to $R_i \ge 0$ called by player one:
\[
\mathcal{W}_{i} = \{ \pi^{-1} (F)  \mid F \in \mathcal{Y}_i \}, \quad 1 \le i \le k.
\]
We have seen that this should give an $R_i$-decomposition of $\mathcal{W}_{i-1}$ over $\mathcal{W}_i$.

If the members of $\mathcal{Y}_k$ are bounded by $T$, the same argument as before shows that for
any subset $F \in \mathcal{Y}_k$ and an element $g_F$ such that $g_F (x_0) \in F$, we have
$g_F^{-1} (F) \subset x_0 [A + B + \ell(T)]$.  So $g_F^{-1} \pi^{-1} (F)$ is contained in the quasi-stabilizer $W = W_{A + 2B + \ell(T)} (x_0)$, and the restriction of $g_F^{-1}$ to $\pi^{-1} (F)$ is an isometry.  

We assume the first player goes on producing numbers $R_{k+1}$, $R_{k+2}$, etc. as part of the game.
Since $W$ has FDC, there is always a winning strategy that can be played entirely (in second player's mind) over $W$ starting with a family $\mathcal{A}_{k+1}$ so that $\{ W \}$ is $R_{k+1}$-decomposable over $\mathcal{A}_{k+1}$ and ending with a bounded family $\mathcal{A}_{k+n}$ for some $n$.  From these auxiliary constructions the second player can produce responses, at every step, to first player's calls as follows:
\[
\mathcal{W}_{i} = \{ g_F \left(g_F^{-1} \pi^{-1} (F) \cap {A} \right) \mid F \in \mathcal{F}_k, A \in \mathcal{A}_i \}, \quad i > k.
\]
The elements $g_F$ act by isometries on $G$, so this gives an $R_i$-decomposition of $\mathcal{W}_{i-1}$ over $\mathcal{W}_i$ for all $i > k$, and if $\mathcal{A}_{k+n}$ is bounded by $U$ then $\mathcal{W}_{k+n}$ is bounded by $U$.

The case of sFDC is entirely similar, with the sequence of numbers $R_i$ being given by player one in advance.
\end{proof}

\SecRef{Applications and discussion}{Apps}

\subsection{}
The following simple corollary to Theorem \refT{Susan} illustrates applications to finitely generated groups that are readily available and require only isometric actions.

\begin{CorRef}{OUR}
Let $\pi \colon G \to H$ be a surjective homomorphism from finitely generated group $G$.  We assume that the groups are given word metrics with respect to finite generating sets and that the kernel $K$ is given the subspace metric.
If $K$ has FAD and $H$ has APC then $G$ has APC.
\end{CorRef}

\begin{proof}
If $S$ is a finite generating set for $G$, $\pi (S)$ can be used as a finite generating set of $H$, and the resulting word metric is known to have APC by quasi-isometry invariance of APC.  The isometric action of $G$ on $H$ is transported from the left action of $G$ on itself.  In this situation, $W_R (e) = K[R]$ from the proof of \cite[Theorem 7]{gBaD:06}.  Since $K$ is a commensurable subspace of all $W_R (e)$, we have $\asdim W_R (e) = \asdim K$.
\end{proof}

\subsection{}

There is a recent effort in extending geometric properties of metric spaces to general coarse properties of coarse structures, motivated by applications in controlled $K$-theory of rings and $C^*$-algebras.  This was done for asymptotic dimension by Grave \cite{bG:06} and for APC and FDC by Bell, Moran, and Nag\'{o}rko \cite{gBdMaN:16}.  The notion of coarse actions is precisely what is needed to formulate group actions on coarse structures, and we expect generalizations of our results to be true when restated for the coarse properties. 

\subsection{}

We could not relax the assumption in Theorem \refT{Susan}.  One can see that the proof relies on the bound for the asymptotic dimension of the quasi-stabilizers known a priori.  This is likely an indication that the more general statement which assumes only that the quasi-stabilizers have FAD is not true.
For much the same reasons, we suspect that the assumptions that $X$ has FAD and quasi-stabilizers have APC do not in general imply that $G$ has APC.

\subsection{}

There exist special positive results of the general type we just dismissed in section 4.3. Guentner \cite{eG:14} points out in 7.2.7 that weak fibered-type conditions for uniformly expansive maps and weak assumptions on the fiber such as simply FAD are sometimes sufficient for extension results.  In his example, the base space $X$ is a simplicial tree, and the action is cofinite. 

We would like to offer a different perspective on the extension problem. 
We will define new \textit{fibred} properties in terms of extensions that are weaker than the absolute analogues but are likely as useful for some purposes. 
The first observation is that there is a natural generalization of coarse quasi-actions. 

\begin{DefRef}{UCAv}
A \textit{(nonuniform) coarse action} of a group $G$ on a metric space $X$ is an assignment of a coarse self-equivalence $f_g \colon X \to X$ to each $g \in G$ so that the following conditions are satisfied:
\begin{enumerate}
\item there is a number $A \ge 0$ such that $d(f_{\id}, \id_X) \le A$ in the sup norm,
\item for each pair of elements $g$ and $h$ in $G$, there is a number $B_{g,h} \ge 0$ such that $d (f_g \circ f_h, f_{gh}) \le B_{g,h}$ in the sup norm.
\end{enumerate}
\end{DefRef}

Just as before, all compositions $f_g \circ f_{g^{-1}}$ are $(A + B_{g,g^{-1}})$-close to the identity.

The reason this is a natural definition is that a coarse action induces a naive $G$-equivariant structure on the bounded $K$-theory spectrum $K(X,R)$.
Each coarse self-equivalence induces a self-equivalence of $K(X,R)$, and the maps close to identity induce genuine identities on $K(X,R)$. So a coarse action on $X$ induces a genuine $G$-action on $K(X,R)$.

On the other hand, we can now make the following definition.

\begin{DefRef}{BMDT}
A finitely generated group $G$ has \textit{fibred property $P_1 \backslash P_2$} if it has a coarse action on a proper metric space $X$ with property $P_1$ so that for all $x \in X$ and all $R \ge 0$ the subsets of $G$ of the form $W_R (x,x_0) = \{ g \in G \mid d (x, gx_0) \le R \}$ have property $P_2$.  
\end{DefRef}

We believe that fibered FAD$\backslash$FAD will be as easy to use in inductive proofs of the integral Novikov Conjecture and the Borel Isomorphism Conjecture in $K$- and $L$-theory as the FAD property itself, cf. \cite{gCbG:04,bG:13}.
So we ask a question of great interest to us.

\begin{QRef}{HJAZX}
How large is the class of groups with fibred property FAD$\backslash$FAD?
It includes groups with FAD.  Does it also include $\Out (F_n)$ with $n > 0$?
Does it include groups with proper isometric actions on a CAT(0) space containing at least one rank-1 element?  Such groups are known to act on quasi-trees.
\end{QRef}

Even when the properties $P_1$ and $P_2$ are both FAD it is unlikely that $G$ has FAD unless one is in the situation of Theorem \refT{Props}.
In the proof of that theorem it was essential that elements of $G$ acted by self-equivalences with the same characteristic function $\ell$.  In the nonuniform case, there is no guarantee that every pull-back $\pi^{-1} (x[R]) = W_R (x,x_0)$ has an isometric embedding in $W_K (x_0)$ for some specific number $K$.  To illustrate this point, we show an example of a well-known infinite dimensional group which has fibred property FAD$\backslash$FAD.

\begin{ExRef}{FADFAD}
The restricted wreath product $\mathbb{Z} \wr \mathbb{Z}$ has FAD$\backslash$FAD.
\end{ExRef}

We will think of the base group as the additive group of the group ring $\mathbb{Z} [\mathbb{Z}]$, so the elements of $\mathbb{Z} \wr \mathbb{Z}$ can be written as pairs $\big( \sum_{x \in \mathbb{Z}} n_x x, k \big)$.  The second factor acts on $\mathbb{Z}$ by translation $k \cdot x = x + k$, so the semidirect product operation is given by $\left( \sum n_x x, k \right) \cdot \left( \sum m_x x, l \right) = \left( \sum (n_x + m_{x-k}) x, k+l \right)$.  
It is known that $\mathbb{Z} \wr \mathbb{Z}$ is generated by two elements.  
We define a nonuniform coarse action by $\mathbb{Z} \wr \mathbb{Z}$ on its cobase $\mathbb{Z}$ by the formula $\left( \sum n_x x, k \right) \cdot t = t + k + n_i / \vert n_i \vert \big( \sum_{x \notin [-t,t]} \vert n_x \vert \big)$, where $i$ stands for the smallest index such that $n_i \ne 0$.  
One easily checks that 
$\asdim W_R (t,0) = {2R+2}$, 
so there is no uniform bound on the dimension of the quasi-stabilizers.  Clearly, $A = 0$.  Using the notation $\lVert \big( \sum n_x x, k \big) \rVert = \lvert k \lvert + \sum \lvert n_x \lvert$, we can choose $B_{g,h} = \lVert g \rVert + \lVert h \rVert$.
Moreover, while the action by each element $\left( \sum n_g g, k \right)$ is eventually the translation by $k$ and so is a coarse equivalence, the function $\ell$ depends linearly on the sum $\sum \vert n_x \vert$, so this is not a uniform coarse quasi-action.

\subsection{}

The notion of fibred property $P_1 \backslash P_2$ does not need to be restricted to groups with geometric actions.
There is the following geometric analogue of Definition \refD{BMDT} in terms of uniformly expansive maps as in \cite{eGrTgY:13} and \cite{gBdM:15}.  

\begin{DefRef}{BMDT2}
A proper metric space $Y$ has \textit{fibred property $P_1 \backslash P_2$} if there is a uniformly expansive map $\pi \colon Y \to X$
where the metric space $X$ has property $P_1$ and, for all $x \in X$ and all $R \ge 0$, the pull-backs $\pi^{-1} (x[R])$ have property $P_2$.  
\end{DefRef}

To point out that this is a useful generalization even when one is interested in geometry of groups, we should recall that the inductive proofs of Novikov and Borel conjectures using controlled algebra are based on the (nonequivariant) Bounded Borel Conjecture which is stated for general metric spaces.    

Now the geometric condition can be naturally iterated.
For example, there is the evident property $P^n = \left( \left( P \backslash P \right)  \backslash P \right) \ldots \backslash P$.

\begin{QRef}{HJAZX}
How large is the class of spaces with property FAD$^n$?  
\end{QRef}

\end{document}